\documentclass[12pt]{amsart}

\newtheorem{thm}{Theorem}
\newtheorem{prop}[thm]{Proposition}

\newtheorem{cor}[thm]{Corollary}

\theoremstyle{remark}
\newtheorem{rem}[thm]{Remark}
\newtheorem{ex}[thm]{Example}

\theoremstyle{definition}

\title{On the simplicial volumes of fiber bundles}
\author{M.~Hoster}
\author{D.~Kotschick}
\address{Mathematisches Institut, Universit\"at M\"unchen,
Theresienstr.~39, 80333 M\"unchen, Germany}
\email{hoster@rz.mathematik.uni-muenchen.de}
\email{dieter@member.ams.org}
\thanks{The second author is grateful to the Institut Mittag-Leffler for
hospitality during the preparation of this paper.}

\begin{document}

\begin{abstract}
We show that surface bundles over surfaces with base and fiber of
genus at least $2$ have non-vanishing simplicial volume.
\end{abstract}

\maketitle


The simplicial volume $\vert\vert M \vert\vert$, introduced by
Gromov~\cite{gromov},
is a homotopy invariant which measures the complexity of the fundamental class
of an oriented manifold $M$. It is determined by the classifying map of
the universal covering, and tends to be non-zero for large manifolds or 
fundamental groups, typically the negatively curved ones.

For products of compact oriented manifolds Gromov~\cite{gromov} proved that
the
simplicial volume is essentially multiplicative. More precisely, there are
universal positive constants $c_n$ depending only on
$n=\text{dim} (M_1\times M_2)$ such that
\begin{equation}\label{gu}
\vert\vert M_1 \times M_2 \vert\vert \leq c_n \vert\vert M_1 \vert \vert \cdot
\vert\vert M_2 \vert\vert \ ,
\end{equation}
and
\begin{equation}\label{gl}
\vert\vert M_1 \times M_2 \vert\vert \geq \vert\vert M_1 \vert \vert \cdot
\vert\vert M_2 \vert\vert \ .
\end{equation}
It is natural to wonder to what extent these inequalities hold for
non-trivial fiber bundles instead of products. In this paper we 
address this question,
initially posed to the second author by A.~Landman.

The analog of the upper bound~\eqref{gu}, which is elementary for products,
fails for non-trivial bundles:

\begin{ex}
There are hyperbolic 3-manifolds $M$ which fiber
over the circle. Like all hyperbolic manifolds, such
$M$ have non-zero simplicial
volume,
whereas $\vert\vert S^1 \vert\vert =0$.
\end{ex}

The lower bound~\eqref{gl}, proved using bounded cohomology, actually
holds for non-trivial bundles if the fiber has dimension 2:

\begin{thm}\label{t:general}
Let $X$ be the total space of a compact oriented
fiber bundle with fiber a compact oriented surface $F$ and  compact
oriented base $B$ (of arbitrary dimension). Then
\begin{equation}\label{l}
\vert\vert X \vert\vert \geq \vert\vert F \vert\vert \cdot
\vert\vert
B \vert\vert \ .
\end{equation}
\end{thm}
\begin{proof}
If $F$ is a sphere or torus, $\vert\vert F\vert\vert =0$, so that there
is nothing to prove. We may thus assume $g(F)\geq 2$.

Let $\pi\colon X\rightarrow B$ be the bundle projection, $e=e(T\pi)\in H^2(X)$
the Euler class of the tangent bundle along the fibers, and $\omega_B\in
H^2(B)$
the fundamental class dual to the orientation class $[B]\in H_2(B)$. 
Let $\pi_*$ denote integration along the fiber. Then
$$
\langle e\cup\pi^*\omega_B , [X]\rangle =\langle\pi_*(e\cup\pi^*\omega_B) , [B]
\rangle =\langle e , [F]\rangle\cdot\langle\omega_B , [B]\rangle = \chi (F) \ .
$$
Denoting by $\vert\vert\cdot\vert\vert_{\infty}$ the Gromov sup norm dual to
the $l^1$ norm used in the definition of the simplicial volume,
cf.~\cite{gromov},
we deduce:
$$
\vert\chi (F) \vert = \vert\vert
e\cup\pi^*\omega_B\vert\vert_{\infty}\cdot\vert\vert
X \vert\vert\leq\vert\vert
e\vert\vert_{\infty}\cdot\vert\vert\omega_B\vert\vert_{\infty}
\cdot\vert\vert X\vert\vert \ ,
$$
or
$$
\vert\vert X\vert\vert\geq\frac{1}{\vert\vert
e\vert\vert_{\infty}}\cdot\frac{1}{\vert\vert
\omega_B\vert\vert_{\infty}}\cdot\vert\chi (F)\vert \ .
$$
As $F$ is assumed to be a hyperbolic Riemann surface, we have 
$\vert\chi (F)\vert=\frac{1}{2}\vert\vert F\vert\vert$, see~\cite{gromov}.
By definition, $\frac{1}{\vert\vert\omega_B\vert\vert_{\infty}}
=\vert\vert B\vert\vert$.
Finally, the unit sphere bundle of $T\pi$ is a flat $S^1$-bundle~\cite{fete},
and so $e$ is bounded. By the Milnor-Wood inequality, 
the exact bound is $\vert\vert e\vert\vert_{\infty}
\leq\frac{1}{2}$ (cf.~\cite{ghys}), so that~\eqref{l} follows.
\end{proof}

\begin{cor}\label{c:1}
Let $X$ be the total space of an oriented surface bundle over a 
surface, with base $F$ and fiber $B$ both of genus $\geq 2$. Then
\begin{equation}\label{ll}
\vert\vert X\vert\vert\geq 4\chi (X) > 0 \ .
\end{equation}
\end{cor}
\begin{proof}
This follows directly from~\eqref{l}, the equality
$\vert\vert\Sigma\vert\vert = -2\chi(\Sigma )$ for closed hyperbolic 
Riemann surfaces $\Sigma$, and the multiplicativity of the 
Euler characteristic in fiber bundles.
\end{proof}

\begin{rem}
Combining~\eqref{ll} with the main result of~\cite{map}, we obtain
$
\vert\vert X\vert\vert \geq 8\vert\sigma (X)\vert$,
where $\sigma$ denotes the signature. However, using $3\sigma (X)=
\langle e^2 , [X]\rangle 
$ and $\vert\vert e\vert\vert_{\infty}\leq\frac{1}{2}$ as in the 
proof of Theorem~\ref{t:general}, we obtain
$$
\vert\vert X\vert\vert\geq 12\vert\sigma (X)\vert \ .
$$
\end{rem}

\begin{cor}
Let $X$ be the total space of a compact oriented fiber bundle with
fiber $F$ and base $B$. If $\text{dim} (X)\leq 4$, then 
\begin{equation}\label{lll}
\vert\vert X\vert\vert\geq\vert\vert F\vert\vert\cdot\vert\vert B\vert\vert \ .
\end{equation}
\end{cor}
\begin{proof}
If $\text{dim}(X) < 4$, then $F$ or $B$ must be a circle, so that the
right-hand-side of~\eqref{lll} vanishes. Similarly, if $\text{dim}(X) = 4$
and either $F$ or $B$ is a circle the right-hand-side
of~\eqref{lll} vanishes and there is nothing to prove. Thus the only
interesting case is that of a surface bundle, for which we appeal
to Theorem~\ref{t:general}.
\end{proof}

It is still an open problem whether there are surface bundles over surfaces
admitting metrics of negative sectional curvature\footnote{Compare
Problem 2.12 (A) in Kirby's list~\cite{kirby}.}. Nevertheless,
Corollary~\ref{c:1} shows that as far as the simplicial volume is 
concerned, surface bundles look like negatively curved manifolds.

\begin{ex}
Kapovich and Leeb~\cite{KL} have given an example of a surface bundle over a
surface with fiber and base both of genus $\geq 2$ which does not admit any
metric of non-positive curvature. By Corollary~\ref{c:1}, the total space has
non-vanishing simplicial volume.
\end{ex}
As far as we know, this is the first example of an aspherical manifold with
non-zero simplicial volume, but with no metric of non-positive curvature. If
one does not insist on asphericity, examples can be constructed using
connected sums, cf.~\cite{gromov}.

Returning now to the upper bound~\eqref{gu} for the simplicial volume, we
have seen already that there is no analog for fibered 3-manifolds. In the
case of fibered 4-manifolds, we do not know if such an inequality holds
in complete generality. However, if we make an additional geometric assumption,
we have:

\begin{prop}\label{p:einstein}
Let $X$ be the total space of a compact oriented fiber bundle with fiber
$F$ and base $B$. If $\text{dim}(X) = 4$, and $X$ admits an Einstein
metric, then
\begin{equation}\label{u}
\vert\vert X \vert\vert \leq c\cdot\vert\vert F \vert\vert \cdot \vert\vert
B \vert\vert \ ,
\end{equation}
for some universal positive constant $c$.
\end{prop}
\begin{proof}
As $X$ is an Einstein manifold, we have Berger's inequality
$\chi (X) \geq 0$, with equality
only if $X$ is flat. Flat manifolds are finitely covered by the 
torus, and so their simplicial volumes vanish and
there is nothing to prove in that case. Thus we can assume $\chi (X)>0$.
This means that neither $F$ nor $B$ can be a circle or a 2-torus. Moreover,
if $F$ is a 2-sphere, then so is $B$, and vice versa. In this case both
sides of~\eqref{u} vanish.

Thus, we only have to consider the case when both $F$ and $B$ are hyperbolic
Riemann surfaces. Then $\vert\vert F\vert\vert\cdot\vert\vert B\vert\vert
=4\chi (F)\cdot\chi (B)=4\chi (X)$, so that we only have to show that the
simplicial volume of $X$ is bounded above by a universal multiple of its
Euler characteristic. But this follows from the Gromov-Hitchin-Thorpe
inequality
$$
\chi (X)\geq\frac{3}{2}\vert\sigma (X)\vert+
\frac{1}{2592\pi^2}\vert\vert X\vert\vert
$$
for Einstein $4$-manifolds proved in~\cite{gromov,GHT}.
\end{proof}

\begin{rem}
The Proposition applies in particular to the surface bundles over surfaces
(with $g(F), g(B) \geq 2$) which admit complex structures, such as
the ones constructed
by Atiyah~\cite{Atiyah} and by Kodaira~\cite{Kodaira}, as they carry
K\"ahler-Einstein metrics due to the results of Aubin and Yau on the
Calabi conjecture.
\end{rem}

The fact that results on surface bundles tend to be stronger when the
total space is complex or admits an Einstein metric is already familiar
from~\cite{map,fibers}.

\bibliographystyle{amsplain}

\begin{thebibliography}{10}

\bibitem{Atiyah}
M.~F.~Atiyah, {\em The signature of fibre bundles}, in
{\sl Global Analysis}, Papers in Honour of K.~Kodaira,
Tokyo University Press 1969.

\bibitem{ghys}
E.~Ghys, {\em Groupes d'homeomorphismes du cercle et cohomologie bornee},
Contemp.~Math. {\bf 58} (1987), 81--106.

\bibitem{gromov}
M.~Gromov, {\em Volume and bounded cohomology},
Publ.~Math.~I.H.E.S. {\bf 56} (1982), 5--99.

\bibitem{KL}
M.~Kapovich and B.~Leeb, {\em Actions of discrete groups on 
nonpositively curved spaces}, Math.~Annalen {\bf 306} (1996), 341--352.

\bibitem{kirby}
R.~Kirby, {\em Problems in low--dimensional topology},
in {\sl Geometric Topology}, ed.~W.~H.~Kazez,
Studies in Advanced Mathematics Vol.~2, Part 2,
American Mathematical Society and International Press 1997.

\bibitem{Kodaira}
K.~Kodaira, {\em A certain type of irregular algebraic surfaces},
J.~Analyse Math. {\bf 19} (1967), 207--215.

\bibitem{GHT}
D.~Kotschick, {\em On the Gromov-Hitchin-Thorpe inequality},
C.~R.~Acad.~Sci.~Paris {\bf 326} (1998), 727--731.

\bibitem{map}
D.~Kotschick, {\em Signatures, monopoles and mapping class groups},
Math.~Research Letters {\bf 5} (1998), 227--234.

\bibitem{fibers}
D.~Kotschick, {\em On regularly fibered complex surfaces}, to appear in
the Kirby Festschrift, a special volume from Geometry and Topology.

\bibitem{fete}
S.~Morita, {\em Characteristic classes of surface bundles and
bounded cohomology}, in {\sl A f$\hat{e}$te of topology}, 
ed.~Y.~Matsumoto et.~al., Academic Press Boston 1988.

\end{thebibliography}

\bigskip

\end{document}